\documentclass[10pt]{article}

\usepackage{multirow}
\usepackage{amsmath}
\usepackage{amsopn}
\usepackage{amsthm}
\usepackage{amsfonts}
\usepackage[T1]{fontenc}
\usepackage[latin1]{inputenc}
\usepackage{graphicx,color}
\usepackage{epsfig}
\DeclareGraphicsExtensions{.eps,.ps}
\usepackage{amssymb}
\usepackage[frenchb,english]{babel}
\usepackage{url}
\usepackage[mathscr]{eucal}
\usepackage{amsbsy}
\usepackage{fancyhdr}
\usepackage{hyperref}

\numberwithin{equation}{section}

\newcounter{rq}
\newcommand\remarque{\refstepcounter{rq}\par\rm{\bf Remark \arabic{rq}.}\quad}
\theoremstyle{plain}
\newtheorem*{Theorem}{Theorem A}

\theoremstyle{plain}
\newtheorem{theoreme}{Theorem}[section]

\newtheorem{definition}[theoreme]{Definition}

\newtheorem{proposition}[theoreme]{Proposition}
\newtheorem{corollaire}[theoreme]{Corollary}
\newtheorem{lemme}[theoreme]{Lemma}
\def\one{\hbox{1\hskip -3pt I}}
\def\rd{\Bbb{R}^d}
\def\sd{S^{d-1}}
\def\qed{\hspace*{2mm} \hfill $\Box$\bigskip}
\def\proof{\par\rm{\it Proof.}\quad}

\def\cok{\hbox{I\hskip -2pt K}}
\def\reip{{\bf P}}

\def\reie{{\bf E}}

\def\cvl1{\stackrel{\mathcal{L}_1}{\longrightarrow}}
\def\cvps{\xrightarrow[n\rightarrow\infty]{a.s.}}

\def\cvpsn{\xrightarrow[N\rightarrow\infty]{a.s.}}

\def\cvas{\xrightarrow[]{a.s.}}

\def\cvdn{\xrightarrow[N\rightarrow\infty]{D}}
\def\cvs{\xrightarrow[n\rightarrow\infty]{}}

\def\egenloi{\stackrel{\mathcal{L}}{=}}

\title{Estimation of parameters of regularly varying distributions on convex cones}
\author{%
\renewcommand{\thefootnote}{\alph{footnote}}
Youri \textsc{Davydov}\,\footnotemark[1]{}\; , %
Shuyan \textsc{Liu}\,\footnotemark[2]{}
}
\date{}
\begin{document}
\maketitle
\renewcommand{\thefootnote}{\alph{footnote}\,}
\footnotetext[1]{Laboratoire Paul Painlevé,
UMR 8524 CNRS Université Lille I, Bât M2, Cité Scientifique,
F-59655 Villeneuve d'Ascq Cedex, France.}

\footnotetext[2]{Institut de statistique, biostatistique et sciences actuarielles,
Université catholique de Louvain,
Voie du Roman Pays, 20, B-1348 Louvain-la-Neuve, Belgium. Research supported by IAP research network grant nr. P6/03 of the Belgian
government (Belgian Science Policy)}

\renewcommand{\thefootnote}{\arabic{footnote}}

\begin{abstract}
The objective of this paper is to extend an estimation method of parameters of the stable distributions in $\rd$ to the regularly varying tails distributions in an arbitrary cone. The consistency and the asymptotic normality of estimators are proved. The sampling method of regrouping is modified to optimize the rate of convergence of estimators.
\end{abstract}

\selectlanguage{frenchb}
\begin{abstract}
Le but de notre travail est de généraliser une méthode d'estimation des paramètres des lois stables dans $\rd$ à des lois à queue régulière dans un cône arbitraire. La consistance et la normalité asymptotique des estimateurs sont prouvées. La méthode d'échantillonnage par paquet est modifiée afin d'optimiser la vitesse de convergence des estimateurs. 
\end{abstract}

\selectlanguage{english}


\newpage


\section{Introduction}
A random $\rd$-valued vector $X$ has a regularly varying tail distribution with characteristic exponent $\alpha>0$ if there exists a finite measure $\sigma$ in the unit sphere $\sd$ such that $\forall B\in\mathcal{B}(\sd)$ with $\sigma(\partial B)=0$,
\begin{equation}\label{vrintro}
\lim\limits_{x\rightarrow\infty}\frac{x^\alpha}{L(x)}\reip\left\{\frac{X}{\|X\|}\in B, \|X\|>x\right\}=\sigma(B),
\end{equation}
where $L$ is a slowly varying function, i.e., $\frac{L(\lambda x)}{L(x)}\rightarrow 1$ as $x\rightarrow\infty, \; \forall \lambda >0$. The measure $\sigma$ is called {\em spectral measure}. Regular variation conditions appear frequently in the studies of the limit theorem for normalised sums of i.i.d. random terms, see e.g. \cite{Rvaceva62} and \cite{Meerschaert01}, and the extreme value theory, see e.g. \cite{Resnick87}. Regular variation is necessary and sufficient conditions for a random $\rd$-valued vector belongs to the domain of attraction of a strictly $\alpha$-stable distribution with $\alpha\in(0,2)$, see e.g. \cite{Araujo80}.

Definition (\ref{vrintro}) of regular variation makes sense in any space, where addition of elements and multiplication by positive scalars are defined, i.e., in any convex cone. Regular variation defined in general metric spaces was studied in \cite{Hult06}. The characterisation of $\alpha$-stable distributions on convex cones and how possible values of the characteristic exponent $\alpha$ relate to the properties of the cone were investigated in \cite{Davydov08}.

We are interested in the problem of estimation of the characteristic exponent $\alpha$ and the spectral measure $\sigma$ of a regularly varying tail distribution on convex cone. By using the relation between the stable distributions and the point processes Davydov and co-workers (see \cite{Davydov99}, \cite{Davydov00a} and \cite{Paulauskas03}) proposed a method to estimate the parameters $\alpha$ and $\sigma$ of the stable distributions in $\rd$. The objective of this work is to extend this method to the regularly varying tails distributions in an arbitrary cone.

Suppose that we have a sample $\xi_1, \xi_2,\ldots, \xi_N$, taken from a regularly varying tail distribution with unknown exponent $\alpha$ and unknown spectral measure $\sigma$. We divide the sample into $n$ groups $G_{m,1},\ldots,G_{m,n}$, each group containing $m$ random elements. In practice, we choose $n=[N^r]$, $r\in (0,1)$, and then $m=[N/n]$, where $[a]$ stands for the integer part of a number $a>0$. As $N$ tends to infinity, we have $nm\sim N$. Let
\begin{equation}\label{defm1}
M^{(1)}_{m,i}=\max\{\|\xi\| \; | \; \xi\in G_{m,i}\}, \ i=1,\ldots,n,
\end{equation}
that is, $M^{(1)}_{m,i}$ denote the largest element in the group $G_{m,i}$. Let $\xi_{m,i}=\xi_j=\xi_{j(m,i)}$ where the index $j(m,i)$ is such that 
\begin{equation}\label{defximi}
\|\xi_{j(m,i)}\|=M^{(1)}_{m,i}.
\end{equation}
We set
\begin{equation}\label{defm2}
M^{(2)}_{m,i}=\max\{\|\xi\| \; | \; \xi\in G_{m,i}\backslash\{\xi_{m,i}\}\}, \ i=1,\ldots,n,
\end{equation}
that is, $M^{(2)}_{m,i}$ denote the second largest element in the same group. Let us denote
\begin{equation}\label{sn}
\varkappa_{m,i}=\frac{M_{m,i}^{(2)}}{M_{m,i}^{(1)}}, \;\;\; \ S_n=\sum_{i=1}^n\varkappa_{m,i}
\end{equation}
and
\begin{equation}\label{defesta}
\hat\alpha_N=\frac{S_n}{n-S_n}.
\end{equation}

The regular variation condition (\ref{vrintro}) implies
\begin{equation}\label{cdqlrd}
\reip\{\|\xi\|>x\}=x^{-\alpha}\tilde{L}(x)+o(x^{-\alpha}) \;\;\;\mbox{as}\;\;\; x\rightarrow\infty
\end{equation}
where $\tilde{L}$ is a slowly varying function. In the sequel we will need the stronger relation : for sufficiently large $x$ and for some $\beta>\alpha$
\begin{equation}\label{cdsdord}
\reip\{\|\xi\|>x\}=C_1x^{-\alpha}+C_2x^{-\beta}+o(x^{-\beta}).
\end{equation}

In the case of $\rd$, under the regular variation assumption and the second-order asymptotic relation (\ref{cdsdord}), the consistency and the asymptotic normality of the estimator $\hat\alpha_N$ were proved firstly for $n=m=[\sqrt{N}]$ in \cite{Davydov00a} and then for more general setting as before in \cite{Paulauskas03}. We resume these results in the following theorem.

\begin{Theorem}\label{thmpaul}(\cite{Paulauskas03})
Let $\xi, \xi_1,\ldots,\xi_N$ be i.i.d. random $\rd$-valued vectors with a distribution satisfying (\ref{cdqlrd}) and let $n=[N^r]$, $r\in (0,1)$, $m=[N/n]$, then
\begin{equation}\label{6.2intro}
\frac{1}{n}S_n\cvpsn\frac{\alpha}{1+\alpha}.
\end{equation}
If the distribution of $\xi$ satisfies (\ref{cdsdord}) with $0<\alpha<\beta\leq\infty$ and we choose
\begin{equation}\label{relmn}
n=N^{2\zeta/(1+2\zeta)-\varepsilon}, \ m=N^{1/(1+2\zeta)+\varepsilon},
\end{equation}
where $\zeta=(\beta-\alpha)/\alpha$ and $\varepsilon\rightarrow 0$ as $N\rightarrow\infty$, then
\begin{equation}\label{6.7}
\frac{\displaystyle\sqrt{n}\left(\frac{1}{n}S_n-\frac{\alpha}{\alpha+1}\right)}{\displaystyle\left(\frac{1}{n}\sum_{i=1}^n\varkappa_{m,i}^2-\left(\frac{1}{n}S_n\right)^2\right)^{1/2}}\Rightarrow
\mathcal{N}(0,1).
\end{equation}
\end{Theorem}

\vspace{0.5cm}
In \cite{Davydov00a} an estimator of normalized spectral measure $\sigma$ was proposed as following. We set 
\begin{equation}\label{deftheta}
\theta_{m,i}=\frac{\xi_{m,i}}{\|\xi_{m,i}\|}, \ i=1,\ldots,n
\end{equation}
where $\xi_{m,i}$ is defined by (\ref{defximi}). Let us denote
\begin{equation}\label{defestms}
\hat{\sigma}_N(\cdot)=\frac{1}{n}\sum_{i=1}^{n}\delta_{\theta_{m,i}}(\cdot).
\end{equation}
Random vectors $\theta_{m,1},\ldots,\theta_{m,n}$ are i.i.d. and it is proved in \cite{Davydov00a} that $\hat{\sigma}_N(\cdot)$ is consistent considering a fixed set $B$, that is, $\forall B\in\mathcal{B}(\sd)$ with $\sigma(\partial B)=0$,
\[\hat{\sigma}_N(B)\cvas\sigma(B).\]
By finding a countable collection of the $\sigma$-continuity sets which is closed under the operation of finite intersection, we obtain $\hat{\sigma}_N\stackrel{a.s.}{\Rightarrow}\sigma \ \mbox{as} \ N\rightarrow\infty$, where $\Rightarrow$ indicate convergence in distribution. The asymptotic normality for $\hat{\sigma}_N(B)$ was proved in \cite{Davydov99}. All these relations were obtained under the assumption that $n=m$. Inspired by the work in \cite{Paulauskas03} we modify the sampling method of regrouping and discuss the convergence rate of the estimator of spectral measure.

Although in the above-mentioned papers the result proved only for a $\rd$-valued sample, it is easy to see that the proposed estimators can be defined in any normed cone and the properties of estimators hold if the relation of the order statistics (6.1) from \cite{Davydov00a} can be verified.

In Section \ref{sect-preliminaires} we summarise several basic definitions related to convex cones and regular variation. The main result of the paper is contained in Section \ref{sect-cons} and \ref{sect-normasymp}. The last section contains the proof of results presented in the previous sections.

\section{Preliminaries}\label{sect-preliminaires}
We summarise several basic definitions related to convex cones and regular variation.  For more detailed introduction we refer the reader to the paper \cite{Davydov08}.

\begin{definition}
An  abelian {\em topological semigroup} is a topological space $\cok$ equipped with a commutative and associative continuous binary operation $+$. It is assumed that $\cok$ \ possesses the neutral element ${\bf e}$ satisfying $x+{\bf e} = x$ for every $x \in \cok$. 
\end{definition}

\begin{definition}\label{defcone}
A {\em convex cone} is an abelian topological semigroup $\cok$ being a metrisable Polish (complete separable) space with a continuous operation $(x, a)\rightarrow ax$ of multiplication by positive scalars for $x\in \cok$ and $a>0$ so that the following conditions are satisfied 
\begin{enumerate}
\item $a(x+y)=ax+ay, \ a>0, \ x,y \in \cok$
\item $a(bx)=(ab)x, \ a,b>0, \ x\in \cok$
\item $1x=x, \ x\in\cok$
\item $a{\bf e}={\bf e}, \ a>0, \ \bf{e} \ \mbox{is the neutral element of \cok}$.
\end{enumerate}

$\cok$  is called a {\em pointed cone} if there is a unique element  {\bf 0} called the {\em origin} such that $ax\rightarrow{\bf 0}$ as $a\downarrow 0$ for any $x \in \cok \ \backslash\{{\bf e}\}$.
\end{definition}




\begin{definition}
A pointed cone \cok \  is said to be a {\em normed cone} if \cok \ is metrisable by a metric $d$ which is homogeneous at the origin, i.e. $d(ax,{\bf 0})=ad(x,{\bf 0})$ for every $a>0$ and $x \in \cok$. The value $\|x\|=d(x,{\bf 0})$ is called the {\em norm} of $x$.
\end{definition}

In the following sections it is assumed that \cok \ is a normed cone. It is obvious that $\|x\|=0$ if and only if $x={\bf 0}$. Furthermore, if ${\bf e}\neq{\bf 0}$, then 4) of Definition \ref{defcone} implies that $\|{\bf e}\|=d({\bf e},{\bf 0})=\infty$. It is therefore essential to allow for $d$ to take infinite values. For instance, if \cok \ is the cone $\overline{\mathbb{R}}_+=[0,\infty]$ with the minimum operation, then the Euclidean distance from any nonempty  $x\in\mathbb{R}_+$ to $\infty$ (being the neutral element) is infinite.


The set
\[S=\{x \; | \; \|x\|=1\}\]
is called the {\em unit sphere}. Note that $S$ is complete with respect to the metric induced by the metric on \cok. The  existence of the origin implies that $\|x\|<\infty$ for all $x\in\cok\ \backslash\{{\bf e}\}$, therefore \cok\  admits a {\em polar decomposition}. This decomposition is realised by the bijection $x\leftrightarrow(\|x\|,x/\|x\|)$ between
\[\cok '=\cok\backslash \{{\bf 0},{\bf e}\}\] 
and $ (0,\infty)\times S$.

In addition to the homogeneity property of the metric $d$, we sometimes require that
\begin{equation}\label{subin}
d(x+h, x)\leq d(h,{\bf 0})=\|h\|, \; x,h\in \cok.
\end{equation}
Then the metric (or the norm) in \cok \ is said to be  {\em sub-invariant}. If \cok \ is a group, then an {\em invariant} (thus also  sub-invariant) metric always exists, i.e. (\ref{subin}) holds with equality sign. This constraint is not trivial, for instance the cone $(\overline{\mathbb{R}}_+, \vee)$ with the maximum operation and Euclidean metric does not satisfy the condition (\ref{subin}).



Before giving the definition of regularity for the random element in a cone \cok \ we recall the definition of a regularly varying function. We say that $L$ is {\em a regular varying function of index $\alpha $ at infinity (respectively at origin)} and we denote $ L \in R_\alpha $ (respectively $ L \in R_\alpha (0+) $) if
\[\frac{L(\lambda x)}{L(x)}\rightarrow x^\alpha,
\ \mbox{as} \  x\rightarrow\infty \; (x\rightarrow 0) \ \mbox{for all} \ \lambda>0.\] 
In particular, if $\alpha=0$ the function $L$ is called {\em slowly varying function}.

\begin{definition}\label{def1}
A random $\cok\,'$-valued element $X$ has a {\em regularly varying tail distribution} if there exists a finite measure $\sigma$ on the unit sphere $S$, a number $\alpha> 0$ and a slowly varying function $L$ such that the convergence (\ref{vrintro}) holds for all $B\in\mathcal{B}(S)$ with $\sigma(\partial B)=0$; here $\|\cdot\|$ is the norm in \cok.
\end{definition}
The measure $\sigma$ is called {\em spectral measure}, and $\alpha$ is called {\em characteristic exponent} or simply {\em tail index}. The fact that $X$ has a regularly varying tail distribution with tail index $\alpha$ and spectral measure $\sigma$ will be noted later by "$X\in \mbox{VR}(\alpha,\sigma)$".



There are various characterizations of the property $X\in \mbox{VR}(\alpha,\sigma)$ (see e.g. \cite{Mikosch03}). We give here an equivalent definition.

\begin{definition}\label{def2}
A random $\cok\,'$-valued element $X\in \mbox{VR}(\alpha,\sigma)$ if there exists a slowly varying function $\tilde L$ such that for all $r>0$ and $B\in\mathcal{B}(S)$ with $\sigma(\partial B)=0$
\begin{equation}\label{regulier2}
\lim_{n\rightarrow\infty}n\reip\left\{\frac{X}{\|X\|}\in B, \|X\|>rb_n\right\}=\sigma(B)r^{-\alpha}, 
\end{equation}
where $b_{n}=n^{1/\alpha}\tilde L(n)$.
\end{definition}

It is well known that in $\rd$ the convergence (\ref{regulier2}) are equivalent to the convergence in distribution of binomial point processes $\beta_n=\sum^{n}_{k=1}\delta_{\xi_k/b_n}$ to a Poisson point process $\pi_{\alpha,\sigma}$ whose intensity measure has a particular form \cite{Resnick87}. This result is generalized to random elements in an abstract cone \cite{Davydov08}. Moreover, in a cone which possesses the sub-invariant norm, the convergence (\ref{regulier2}) with $\alpha\in (0, 1)$ implies that $X$ belongs to the domain of attraction of a strictly $\alpha$-stable distribution (see Th. 4.7 \cite{Davydov08}).

Typical examples of cones that fulfil our requirements are Banach spaces or convex subcones in Banach spaces; the family of compact (or convex compact) subsets of a Banach space with Minkowski addition, see e.g. \cite{Davydov00a} \cite{Gine85} \cite{Gine82}; the family of compact sets in $\rd$ with the union operation, see e.g. \cite{Molchanov05}; the family of all finite measures with the conventional addition operation and multiplication by numbers, see e.g. \cite{Daley03} \cite{Rachev91}. Another typical example is the set $\mathbb{R}_+=[0,\infty)$ with the maximum operation, i.e. $x+y=x\vee y=\max(x,y)$. Max-stable distributions appear in this case. More information on these and other examples can be found in \cite{Davydov08}.
\section{The consistency of estimators}\label{sect-cons}
We assume that $\xi_1, \xi_2,\ldots, \xi_N$ are i.i.d. random $\cok$- valued elements with a regularly varying tail distribution. For our purposes we suppose that the spectral measure is normalized, i.e., $\sigma(S)=1$. Our aim is to estimate $\alpha$ and $\sigma$ from the sample. The estimators $\hat\alpha_N$ and $\hat\sigma_N(\cdot)$ are defined by (\ref{defesta}) and (\ref{defestms}). To establish consistency of these estimators we need two auxiliary results below.

Let $X$ be a $\cok\,'$-valued random element, we denote $G(x)=\reip\{\|X\|>x\}$. Let $Y_1, Y_2, \ldots$ be the random variables i.i.d. with distribution function $1-G$ and $Y_{n,1}, Y_{n,2}, \cdots , Y_{n,n}$, $Y_{n,1}\geq Y_{n,2}\geq \cdots \geq Y_{n,n}$, the corresponding order statistics. 
\begin{lemme}\label{bny}
If $X$ satisfies the condition of regular variation (\ref{regulier2}), then the vector $b_n^{-1}(Y_{n,1},\ldots,Y_{n,n},0,0,\ldots)$ converge in distribution in $\mathbb{R}^\infty$ to\\ 
$\sigma(S)^{1/\alpha}(\Gamma_1^{-1/\alpha},\Gamma_2^{-1/\alpha},\ldots)$, where $\Gamma_i=\sum^i\limits_{j=1}\lambda_j$ and $\lambda_1, \lambda_2,\ldots$ are i.i.d. random variables with a standard exponential distribution, i.e. $\reie (\lambda_{i})=1$.
\end{lemme}
The proof is given in Section \ref{sect-proof}.   

The second result is a variant of the strong law of large numbers for a triangular array.
\begin{proposition}\label{lemgnc2}
Let $\{X_{m,i}, 1\leq i\leq n\}$ be i.i.d. real random variables for each $m$. Suppose that the indices $n$ and $m$ satisfy the following relations
\begin{equation}\label{relmnlemgnc2}
n\sim N^r, \;\;\; m\sim N^{1-r} \;\;\; \mbox{as} \;\;\; N\rightarrow\infty
\end{equation}
where $0<r<1$ is a constant and $N\in\Bbb{N}$. If there exists a real number $k>\frac{2}{r}$ and a constant $M>0$ such that $\reie|X_{m,1}|^{k}\leq M<\infty$, then
\begin{equation}\label{lgnc2}
\frac{1}{n}\sum_{i=1}^{n}X_{m,i}-\reie X_{m,1}\xrightarrow[N\rightarrow\infty]{a.s.} 0.
\end{equation}
\end{proposition}
The proof is given in Section \ref{sect-proof}.
\vspace{5mm}

We consider now the estimator of characteristic exponent $\hat\alpha_N$.
\begin{theoreme}\label{sncvpsc}
Let $\xi, \xi_1,\ldots,\xi_N$ be i.i.d. $\cok\, '$-valued random elements and $\xi\in \mbox{VR}(\alpha,\sigma)$. If $S_n$ is defined by (\ref{sn}) with $n\sim N^r$, $0<r<1$, then
\begin{equation}\label{6.2}
\frac{1}{n}S_n\cvpsn\frac{\alpha}{1+\alpha}.
\end{equation}
\end{theoreme}

\remarque This means that the quantity $\frac{\displaystyle S_n}{\displaystyle n-S_n}$ gives a consistent estimator of the parameter $\alpha$.
\vspace{0.5cm}

\proof
By Lemma \ref{bny}, it follows that for all $i$
\begin{equation}\label{6.1}
\left(\frac{M^{(1)}_{m,i}}{b_m},\frac{M^{(2)}_{m,i}}{b_m}\right)\Rightarrow c(\Gamma_1^{-1/\alpha},\Gamma_2^{-1/\alpha}) \;\;\;  \mbox{as} \;\;\;  m\rightarrow\infty.
\end{equation}
Therefore for all $i, 1\leq i \leq n$ we have
\[\varkappa_{m,i}\Rightarrow\left(\frac{\lambda_1}{\lambda_1+\lambda_2}\right)^{1/\alpha} \; \mbox{as}\;\;\; m\rightarrow\infty.\]
Since $0\leq\varkappa_{m,i}\leq 1$, for any integer $k$
\begin{equation}\label{6.3}
\reie\varkappa_{m,i}^k\xrightarrow [m \rightarrow
\infty]{}\reie\left(\frac{\lambda_1}{\lambda_1+\lambda_2}\right)^{k/\alpha}.
\end{equation}
Simple calculations show that  
\begin{equation}\label{6.4}
\reie\left(\frac{\lambda_1}{\lambda_1+\lambda_2}\right)^{1/\alpha}=\frac{\alpha}{\alpha+1}.
\end{equation}
The random variables $\varkappa_{m,1},\ldots,\varkappa_{m,n}$ are i.i.d. and $0\leq\varkappa_{m,i}\leq 1.$ By applying Proposition \ref{lemgnc2}, we have
\begin{equation}\label{sncve}
\frac{1}{n}S_n-\reie\varkappa_{m,1}\cvpsn 0.
\end{equation}
 Combining this with (\ref{6.3}) (\ref{6.4}) proves (\ref{6.2}), completing the proof of the theorem. \qed
 \vspace{0.5cm}

Next we consider the estimator of the spectral measure defined by (\ref{defestms}). Note that the random elements $\theta_{m,1},\ldots,\theta_{m,n}$ on which this estimator is defined are i.i.d. in $S$. The following lemma shows the asymptotic property for each $\theta_{m,i}$, $i=1,\ldots,n$.
\begin{lemme}\label{cvdetheta}
Let $\xi_1,\ldots,\xi_N$ be i.i.d. $\cok\, '$-valued random elements such that the condition (\ref{regulier2}) is satisfied. If $\theta_{m,i}$ is defined by (\ref{deftheta}), then
\begin{equation}\label{thetacv}
\theta_{m,i}\Rightarrow\sigma \ \mbox{as} \ m\rightarrow\infty
\end{equation}
for each $i$.
\end{lemme}
\proof For all Borel set $B$ in unite sphere $S$ such that $\sigma(\partial B)=0$, we have
\begin{eqnarray*}
\reip\{\theta_{m,i}\in B\}&=&\reip\{{\xi_{m,1}}/{\|\xi_{m,1}\|}\in B\}\\
&=&\sum_{k=1}^m\reip\{{\xi_{m,1}}/{\|\xi_{m,1}\|}\in B, \xi_{m,1}={\xi_k}\}\\
&=&m\reip\{{\xi_{m}}/{\|\xi_{m}\|}\in B, \xi_{m,1}=\xi_m\}\\
&=&m\reip\{{\xi_m}/{\|\xi_m\|}\in B, {\|\xi_m\|}\geq{\|\xi_k\|}, \forall k=1,\ldots,m-1\}\\
&=&m\reip\{{\xi_m}/{\|\xi_m\|}\in B, {\|\xi_m\|}\geq b_m\tau_{m-1}\}\\
&=&\int m\reip\{{\xi_m}/{\|\xi_m\|}\in B, \|\xi_m\|\geq b_mx\}\reip_{\tau_{m-1}}(dx).
\end{eqnarray*}
where $\reip_{\tau_{m-1}}$ is the distribution of $\tau_{m-1}=\max\limits_{1\leq k\leq m-1}(\|\xi_k\|b_m^{-1})$.

By (\ref{regulier2}) and Lemma \ref{bny}, the last term converges to
\[\int \sigma(B)x^{-\alpha}\reip_{\Gamma_1^{-1/\alpha}}(dx)=\sigma(B)\reie(\Gamma_1^{-1/\alpha})^{-\alpha}=\sigma(B).\]
\qed

Therefore for each Borel set $B$ in unit sphere $S$ such that $\sigma(\partial B)=0$, we have
\[\one_B(\theta_{m,i})\Rightarrow\one_B(\eta) \; \mbox{as} \; m\rightarrow\infty,\]
where $\eta$ is a random element with distribution $\sigma$. This yields
\begin{equation}\label{6.5}
\reie\one_B(\theta_{m,i})\xrightarrow[m\rightarrow\infty]{}\sigma (B).
\end{equation}

If there exists a constant $r>0$ such that $n\sim N^r$, applying Proposition \ref{lemgnc2} for the triangular array $\{\one_B(\theta_{m,1}),\ldots,\one_B(\theta_{m,n})\}$, we have for each fixed set $B\in\mathcal{B}(S)$ with $\sigma(\partial B)=0$,
\begin{equation}\label{6.6}
\frac{1}{n}\sum_{i=1}^{n}\one_B(\theta_{m,i})-\reie\one_B(\theta_{m,1})\cvpsn 0.
\end{equation}
Together (\ref{6.5}) and (\ref{6.6}) we have the following result.

\begin{theoreme}\label{connuc}
Let $\xi, \xi_1,\ldots,\xi_N$ be i.i.d. $\cok\, '$-valued random elements and $\xi\in \mbox{VR}(\alpha,\sigma)$. If $\hat\sigma_N(\cdot)$ is defined by (\ref{defestms}) with $n\sim N^r$, $0<r<1$, then $\forall B\in\mathcal{B}(S)$ with $\sigma(\partial B)=0$,
\begin{equation}\label{defmuc}
\hat{\sigma}_N(B)=\frac{1}{n}\sum_{i=1}^{n}\one_B(\theta_{m,i})\cvpsn\sigma(B).
\end{equation}
\end{theoreme}

The result is for a fixed set. A stronger convergence can be proved by an immediate application of the following proposition. 
\begin{proposition}\label{lempscv}
Let $(S,\mathcal{S})$ be a complete separable metric space. Let $\{\sigma_n\}$ be a sequence of random probability measure in $S$. If $\sigma$ is a probability measure on $(S,\mathcal{S})$ such that for each set $B\in\mathcal{B}(S)$ with $\sigma(\partial B)=0$ we have the convergence $\sigma_n(B)\stackrel{a.s.}{\rightarrow}\sigma(B),$ then
\[\sigma_n\stackrel{a.s.}{\Rightarrow}\sigma \ \mbox{as} \ n\rightarrow\infty.\]
\end{proposition}
The proof is given in Section \ref{sect-proof}

\begin{corollaire}\label{cvsigmahat}
Under the same assumption of Theorem \ref{connuc}, we have
\[\hat{\sigma}_N\stackrel{a.s.}{\Rightarrow}\sigma \ \mbox{as} \ N\rightarrow\infty.\]
\end{corollaire}

\section{The asymptotic normality of estimators}\label{sect-normasymp}
In this section we consider the asymptotic normality of the estimators $\hat{\alpha}_N$ and $\hat{{\sigma}}_N$. 

Let us denote 
\[a_m=\reie\varkappa_{m,1}=\frac{\alpha}{\alpha+1}+r_m,\]
where $r_m\rightarrow 0$ as $m\rightarrow\infty$. Using this notation, we can write
\[\sqrt{n}\left(\displaystyle\frac{1}{n}S_n-\frac{\alpha}{\alpha+1}\right)=\frac{1}{\sqrt{n}}\sum_{i=1}^n(\varkappa_{m,i}-a_m)+\sqrt{n}r_m.\]
It is easy to see that (\ref{6.7}) follows from the following three relations :
\begin{equation}\label{6.8}
\frac{1}{\sqrt{n}}\sum_{i=1}^n(\varkappa_{m,i}-a_m)\Rightarrow
\mathcal{N}(0,\sigma^2),
\end{equation}

\begin{equation}\label{6.9}
\sqrt{n}r_m\rightarrow 0,
\end{equation}

\begin{equation}\label{6.10}
\frac{1}{n}\sum_{i=1}^n\varkappa_{m,i}^2-\left(\frac{1}{n}S_n\right)^2\xrightarrow[N\rightarrow\infty]{P}\sigma^2.
\end{equation}

Denoting $\sigma_m^2:=\reie(\varkappa_{m,i}-a_m)^2$, we have
\[\sigma_m^2\rightarrow\sigma^2=\alpha((\alpha+1)^2(\alpha+2))^{-1}.\]
The relation is a consequence of (\ref{6.3}). Since the i.i.d. random variables $\{\varkappa_{m,i}-a_m, 1\leq i\leq n\}$ are uniformly bounded, we assume $|\varkappa_{m,i}-a_m|\leq M<\infty$. Considering the condition of Lindeberg : for all $\varepsilon>0$,
\begin{eqnarray*}
&&\lim_{n\rightarrow\infty}\frac{1}{n\sigma_m^2}\sum_{i=1}^n\int_{\{|\varkappa_{m,i}-a_m|>\varepsilon\sqrt{n}\sigma_m\}}(\varkappa_{m,i}-a_m)^2d\reip\nonumber\\
&\leq&\lim_{n\rightarrow\infty}\frac{M^2}{\sigma_m^2}\reip\{|\varkappa_{m,1}-a_m|>\varepsilon\sqrt{n}\sigma_m\}\nonumber\\
&=&0,
\end{eqnarray*}
the relation (\ref{6.8}) follows from the central limit theorem applied to triangular array $\{\varkappa_{m,i}-a_m, 1\leq i\leq n\}$. The relation (\ref{6.10}) follows from Proposition \ref{lemgnc2} applied to triangular array $\{\varkappa_{m,i}^2, 1\leq i\leq n\}$ and Theorem \ref{sncvpsc}. Therefore it remains to prove (\ref{6.9}). The same proof of Theorem 1 in \cite{Paulauskas03} gives the following result.

\begin{theoreme}\label{asyma}
Let $\xi, \xi_1,\ldots,\xi_N$ be i.i.d. $\cok\, '$-valued random elements and $\xi\in \mbox{VR}(\alpha,\sigma)$. If the distribution of $\xi$ satisfies the relation (\ref{cdsdord}) with $0<\alpha<\beta\leq\infty$ and we choose $n$ and $m$ verifying the relation defined by (\ref{relmn}), then we have (\ref{6.7}).
\end{theoreme}

Using Delta-method we obtain the result in a more convenient form for the statistics.
\begin{corollaire}\label{cornormac}
Under the same assumption of Theorem \ref{asyma} we have
\begin{equation}\label{cor6.7}
\frac{\displaystyle\sqrt{n}\left(\hat{\alpha}_N-\alpha\right)}{(\alpha+1)^2\displaystyle\left(\frac{1}{n}\sum_{i=1}^n\varkappa_{m,i}^2-\left(\frac{1}{n}S_n\right)^2\right)^{1/2}}\Rightarrow
\mathcal{N}(0,1).
\end{equation}
\end{corollaire}
\remarque
Simple calculations show that
\[\frac{1}{n}\sum_{i=1}^n\varkappa_{m,i}^2-\left(\frac{1}{n}S_n\right)^2\xrightarrow[N\rightarrow\infty]{P}\alpha((\alpha+1)^2(\alpha+2))^{-1}.\]
Therefore the asymptotic variance of estimator $\hat\alpha_N$ is $\frac{\alpha(\alpha+1)^2}{\alpha+2}$. This means that the confidence intervals of $\hat{\alpha}_N$ increase with $\alpha$.

\vspace{0.5cm}

Before considering the asymptotic normality of the estimator of spectral measure, we present a {\em strong second-asymptotic relation} : $\forall B\in\mathcal{B}(S)$ with $\sigma(\partial B)=0$,
\begin{equation}\label{2ndm}
\reip\left\{\frac{\xi}{\|\xi\|}\in B, \|\xi\|>x\right\}=\sigma(B)x^{-\alpha}+Cx^{-\beta}+o(x^{-\beta}) \;  \mbox{as} \;
x\rightarrow\infty,
\end{equation}
where $\beta >\alpha >0$. Note that this condition implies the conditions (\ref{vrintro}) and (\ref{cdsdord}).

Let us denote $\sigma(B)=b, \  \one_B(\theta_{m,i})=\eta_{m,i}, i=1,2,\ldots,n$. Then
\[Z_n=\hat{\sigma}_N(B)-\sigma(B)=\frac{1}{n}\sum_{i=1}^n(\eta_{m,i}-\sigma(B)),\]
\[\sqrt{n}Z_n=\frac{1}{\sqrt{n}}\sum_{i=1}^n(\eta_{m,i}-\sigma(B))=\frac{1}{\sqrt{n}}\sum_{i=1}^n(\eta_{m,i}-\reie\eta_{m,1})+\sqrt{n}(\reie\eta_{m,1}-\sigma(B)),\]
\[U_{n}=\frac{1}{\sqrt{n}}\sum_{i=1}^n(\eta_{m,i}-\reie\eta_{m,1}), \   r_{m}=\reie\eta_{m,1}-\sigma(B).\]

We set
\[T_N=\displaystyle\frac{\sqrt{n}Z_n}{\left(\displaystyle\frac{1}{n}\sum_{i=1}^n\eta_{m,i}^2-\left(\displaystyle\frac{1}{n}\sum_{i=1}^n\eta_{m,i}\right)^2\right)^{1/2}},\]
then the asymptotic property of $\hat{\sigma}_N$ can be described as following.

\begin{theoreme}\label{normademc}
Let $\xi, \xi_1,\ldots,\xi_N$ be i.i.d. $\cok\, '$-valued random elements and the distribution of $\xi$ satisfies the condition (\ref{2ndm}). If we choose
\[n=N^{2\zeta/(1+2\zeta)-\varepsilon},\;\; m=N^{1/(1+2\zeta)+\varepsilon},\]
where $\zeta=\min(\frac{\beta-\alpha}{\alpha}, 1)$ and $\varepsilon$ is an arbitrarily small positive constant, then
\begin{equation}\label{29}
T_N\Rightarrow \mathcal{N}(0,1).
\end{equation}
\end{theoreme}

\remarque 
We can get also asymptotic normality for $\sqrt{n}Z_n$, but the variance of the limit normal law is $\sigma(B)(1-\sigma(B))$ which we are estimating.

\remarque
Fristedt, \cite{Fristedt72}, proved an asymptotic expansion for the distribution of the norm of a strictly $\alpha$-stable random vector in $\rd$
\begin{equation}\label{frist}
G(x)=c_1x^{-\alpha}+c_2x^{-2\alpha}+O(x^{-3\alpha}), \ \mbox{as}\; x\rightarrow\infty.
\end{equation}
That means $\beta=2\alpha$. In this case the optimal choice of $n$ and $m$ is $n=N^{2/3-\varepsilon}, m=N^{1/3+\varepsilon}$, hence the rate of convergence of $\hat{\sigma}_N(B)$ in $\mathcal{L}_{1}$ is close to $N^{1/3}$.
\vspace{0.5cm}

The proof is given in Section \ref{sect-proof}.

\section{Proofs}\label{sect-proof}

\par{\it Preliminary remarks}\\

Let $X$ be $\cok'$-valued random element satisfying the condition of regular variation (\ref{regulier2}). We denote $Y=\|X\|$ and $G(x)=\reip\{Y>x\}$. Then
\begin{equation}\label{nG}
nG(b_nx)\rightarrow \sigma(S)x^{-\alpha},\;  \mbox{as} \; n\rightarrow \infty, \; \mbox{for all} \; x>0.
\end{equation}
For positive fixed $x$, we choose $n$ the smallest integer such that $b_{n+1}>x$. Then $b_n\leq x<b_{n+1}$ and for a non-creasing function $G$ we have
\[\frac{G(\lambda b_{n+1})}{G(b_n)}\leq\frac{G(\lambda x)}{G(x)}\leq\frac{G(\lambda b_{n})}{G(b_{n+1})},\; \; \mbox{for all}\; \lambda>0.\]
By (\ref{nG}) we have $nG(b_n)\rightarrow\sigma(S)$, then
\[\frac{G(\lambda x)}{G(x)}\rightarrow\lambda^{-\alpha}, \; \mbox{as} \; x\rightarrow\infty\;  \mbox{for all}\; \lambda>0.\]
We deduce that $G\in R_{-\alpha}$, which allow us write the following equivalence
\begin{equation}\label{G}
G(x)\sim x^{-\alpha}L(x),
\end{equation}
where $L(x)$ is a slowly varying function.

We recall a well known result on the asymptotic inverse of a regular varying function.
\begin{theoreme}\label{variation} (\cite{Bingham87} Th. 1.5.12)
Let $f\in R_\alpha$ with $\alpha >0$, then $\exists g(x)\in R_{1/\alpha}$ such that the following relation holds 
\begin{equation}\label{equiva}
f(g(x))\sim g(f(x))\sim x \; \mbox{as} \; x\rightarrow\infty.
\end{equation}
Here $g$ ({\em the asymptotic inverse} of $f$) is defined uniquely up to asymptotic equivalence, and a version of $g$ is 
\[f^{\leftarrow}(x)=\inf\{y: f(y)\leq x\}.\]
\end{theoreme}

We denote
\begin{equation}\label{eqvauxthm1}
f(x)=\frac{1}{G(x)}\sim x^\alpha\frac{1}{L(x)},
\end{equation}
then $f(x)\in R_\alpha$ with $\alpha>0$. By applying the previous theorem, we obtain the inverse $g(x)$ of $f(x)$ in the following form,
\[g(x)=x^{1/\alpha}L^\sharp(x)\]
where the slowly varying function $L^\sharp$ verifies the following relation 
\begin{equation}\label{l1}
L(x)^{-1/\alpha}L^\sharp(f(x))\rightarrow 1,
\end{equation}
and
\begin{equation}\label{l3}
L(g(x))^{-1}L^\sharp(x)^\alpha\rightarrow 1,\ x\rightarrow\infty .
\end{equation}
By (\ref{equiva}) and (\ref{eqvauxthm1}) we have
\begin{equation}\label{equideg}
G\left(g\left(\frac{1}{x}\right)\right)\sim x, \ x\rightarrow 0.
\end{equation}
Defining the generalized inverse
\[G^{-1}(x):=\inf\{y: G(y)< x\},\] 
we can prove that
\begin{equation}\label{equideg2}
G(G^{-1}(x))\sim x, \ x\rightarrow 0.
\end{equation}
For this we choose $\lambda>1$, $A>1$, $\delta\in(0,\infty)$, then by the theorem of Potter (Th. 1.5.6 \cite{Bingham87} page 25) there exists $u_0$ such that
\[A^{-1}\lambda^{-\alpha-\delta}G(v)\leq G(u)\leq A\lambda^{\alpha+\delta}G(v), \; \forall v\in[\lambda^{-1}u, \lambda u], \;\; u\geq u_0.\]
We take $x$ small enough such that $G^{-1}(x)\geq u_0$, then by the definition of $G^{-1}$ there exists $y\in[\lambda^{-1}G^{-1}(x),G^{-1}(x)]$ such that $G(y)\geq x$, and there exists $y'\in[G^{-1}(x),\lambda G^{-1}(x)]$ such that $G(y')< x$. Taking $G^{-1}(x)$ for $u$, $y$ and $y'$ for $v$, we get
\[A^{-1}\lambda^{-\alpha -\delta}G(y)\leq G(G^{-1}(x))\leq A\lambda^{\alpha+\delta}G(y').\]
Hence $\limsup$ and $\liminf$ of $G(G^{-1}(x))/x$ are between $A\lambda^{\alpha+\delta}$ and $A^{-1}\lambda^{-\alpha-\delta}$  as $x\rightarrow \infty$. Taking $A$, $\lambda \downarrow 1$, we have $G(G^{-1}(x))/x\rightarrow 1$. 

The relations (\ref{equideg}) and (\ref{equideg2}) give immediately
\[G^{-1}(x)\sim g(1/x), \ x\rightarrow 0.\]  
Thus we have the equivalent expression of the inverse of $G(x)$:
\begin{equation}\label{G-1}
G^{-1}(x)\sim x^{-1/\alpha}L^\sharp(1/x)\in R_{-1/\alpha}(0+).
\end{equation}

\begin{lemme}
Let $X$ be a $\cok\,'$-valued random element, $G(x)=\reip\{\|X\|>x\}$. If the distribution of $X$ satisfies the condition (\ref{regulier2}), then for each
$i=1,2,\ldots,$
\begin{equation}\label{inverg}  
b_n^{-1}G^{-1}\left(\frac{\Gamma_i}{\Gamma_{n+1}}\right)\cvs\sigma(S)^{1/\alpha}\Gamma_i^{-1/\alpha} \ \mbox{with probability 1.}
\end{equation}
\end{lemme}

\proof We recall (\ref{nG}) $nG(b_nx)\rightarrow\sigma(S)x^{-\alpha}$ which implies 
\[G(x_n)\sim\frac{\sigma(S)}{n}\left(\frac{b_n}{x_n}\right)^{\alpha}, \; n\rightarrow\infty\] 
where $x_n=b_n x$, $x_n\rightarrow\infty$, as $n\rightarrow\infty$. By replacing the left term in the previous formula by (\ref{G}), we get an equivalent expression of $b_n$ in terms of $L(x)$:
\begin{equation}\label{bnva}
b_n\sim \left(\frac{nL(x_n)}{\sigma(S)}\right)^{1/\alpha}, \
n\rightarrow\infty.
\end{equation}

Considering (\ref{G-1}), we have an equivalent expression with probability 1 for each $i$,
\begin{equation}\label{g-1}
G^{-1}\left(\frac{\Gamma_i}{\Gamma_{n+1}}\right)\sim\left(\frac{\Gamma_i}{\Gamma_{n+1}}\right)^{-1/\alpha}L^\sharp\left(\frac{\Gamma_{n+1}}{\Gamma_{i}}\right),
\ n\rightarrow\infty
\end{equation}
where $L^\sharp$ satisfies (\ref{l1}) which means
\begin{equation}\label{l2}
L(x_n)^{-1/\alpha}L^\sharp\left(\frac{\Gamma_{n+1}}{\Gamma_{i}}\right)\rightarrow
1,\ n\rightarrow\infty .
\end{equation}
Collecting (\ref{bnva})-(\ref{l2}) we deduce that with probability $1$
\[b_n^{-1}G^{-1}\left(\frac{\Gamma_i}{\Gamma_{n+1}}\right)\sim\sigma(S)^{1/\alpha}\left(\frac{\Gamma_{n+1}}{n\Gamma_i}\right)^{1/\alpha}L(x_n)^{-1/\alpha}L^\sharp\left(\frac{\Gamma_{n+1}}{\Gamma_i}\right)\sim \left(\frac{\sigma(S)}{\Gamma_i}\right)^{1/\alpha},\]
as $n\rightarrow\infty$, since $\Gamma_{n+1}/n\cvps 1$; the lemma is proved.
\qed

\par{\it Proof of Lemma \ref{bny}:}\quad
It is well known that (see, e.g. \cite{Breiman68} Section 13.6)
\begin{equation}\label{yordre}
(Y_{n,1},\ldots,Y_{n,n})\egenloi \left(G^{-1}\left(\frac{\Gamma_1}{\Gamma_{n+1}}\right),\ldots,G^{-1}\left(\frac{\Gamma_n}{\Gamma_{n+1}}\right)\right).
\end{equation}
The lemma follows from (\ref{inverg}) and (\ref{yordre}).\qed

\par{\it Proof of Proposition \ref{lemgnc2}:}\quad
Denote $Y_{m,i}=X_{m,i}-\reie X_{m,1}$, then the random variables $\{Y_{m,i}, 1\leq i\leq n\}$ are centered and i.i.d.. We have
\[\reie|Y_{m,1}|^{k}\leq \reie (|X_{m,1}|+|\reie X_{m,1}|)^{k}\leq\reie(2^{k-1}(|X_{m,1}|^{k}+|\reie X_{m,1}|^{k}))\leq 2^{k}\reie|X_{m,1}|^{k}.\]
It is well  known (see \cite{Rosenthal70}) that for $k\geq 2$ we have
\[\reie\left|\sum_{i=1}^nY_{m,i}\right|^{k}\leq c(k)n^{k/2}\reie|Y_{m,1}|^{k},\] 
where $c(k)$ is a positive constant depending only on $k$. We deduce from the condition (\ref{relmnlemgnc2}) that there exists a constant $C>0$ such that $n\geq CN^r$. Hence for all $\varepsilon >0$ we have
\[\reip\left\{\left|\frac{1}{n}\sum_{i=1}^nY_{m,i}\right|>\varepsilon\right\}\leq\frac{\reie\left|{\displaystyle \sum_{i=1}^nY_{m,i}}\right|^{k}}{n^{k}\varepsilon^{k}}\leq \frac{2^{k}c(k)\reie|X_{m,1}|^{k}}{C^{\frac{k}{2}}N^{\frac{kr}{2}}\varepsilon^{k}}=\frac{c_0}{N^{\frac{kr}{2}}\varepsilon^{k}},\]
where $c_0=2^{k}c(k)M/C^{\frac{k}{2}}$. Since $\frac{kr}{2}>1$, we can find a small enough positive number $\varepsilon'$ such that $\frac{kr}{2}-\varepsilon '>1$. Taking $\varepsilon=\varepsilon_N=N^{-\frac{\varepsilon '}{k}}$ and applying the Borel-Cantelli lemma, we get that with probability $1$ and for $N$ large enough
\[\left|\frac{1}{n}\sum_{i=1}^nX_{m,i}-\reie X_{m,1}\right|\leq N^{-\frac{\varepsilon '}{k}},\]
the proposition is proved. 
\qed

\par{\it Proof of Proposition \ref{lempscv}:}\quad
We denote the collection of all $\sigma$-continuity sets by
\[\mathcal{D}_\sigma=\{B~|~B\in\mathcal{B}(S),\;\sigma(\partial B)=0\}.\] 
Since $S$ est separable, there exists a countable dense set in $S$, denoted by
\[W=\{x_1,x_2,\ldots\},\ x_i\in S, \ i=1,2,\ldots.\]
We denote the open ball with centre $x_i$ in $W$ and radius $r$ by 
\[V(x_i,r)=\{x~|~ x\in S, \|x-x_i\|<r\}.\]
Since for each $x_i\in W$ the boundaries $\partial\{V(x_i,r)\}\subset\{x~|~\|x-x_i\|=r\}$ are disjoints for different $r$, at most a countable number of them can have positive $\sigma$-measure. Therefore, there exists a sequence of positive numbers $r_k^i \downarrow 0$ as $k\rightarrow\infty$ for each $x_i$ such that
\[\mathbb{L}_{i}=\{V(x_i,r_k^i), k=1,2,\ldots\}\subset\mathcal{D}_\sigma.\]
The collection $\mathbb{L}=\bigcup\limits_{x_i\in W}\mathbb{L}_{i}$ is countable. It is clear that for each $x_i\in W$, the collection $\mathbb{L}_{i}$ is a local base at point $x_i$ for la topology $\mathcal{S}$. Since $W$ is dense in $S$, $\mathbb{L}$ is a base of $\mathcal{S}$. The $\sigma$-algebra generated by $\mathbb{L}$, denoted by $\sigma(\mathbb{L})$, is the Borel-field $\mathcal{B}(S)$.

Now we expand $\mathbb{L}$ by adding the finite intersections of members of $\mathbb{L}$, we denote
\[\mathcal{L}=\mathbb{L}\cup\left.\left\{\bigcap_{i\in I}V_{i}~\right|~V_{i}\in\mathbb{L}, I\subset\mathbb{N}, \;\mbox{card}(I)<\infty\right\}.\]
It is clear that $\mathcal{L}$ is still countable and $\sigma(\mathcal{L})=\mathcal{B}(S)$, moreover $\mathcal{L}\subset\mathcal{D}_\sigma$. Since $\sigma_n(B)\stackrel{a.s.}{\rightarrow}\sigma(B)$ for all $B\in\mathcal{B}(S)$ and $\sigma(\partial B)=0$, then $\forall V\in\mathcal{L}$, $\exists\Lambda_V\subset\Omega$ and $\reip(\Lambda_V)=0$, such that $\forall \omega\in\Lambda_V^\complement$ we have
\begin{equation}\label{cvsigma}
\sigma_n(\omega,V)\rightarrow\sigma(V).
\end{equation}
If we denote $\Lambda=\bigcup\limits_{V\in\mathcal{L}}\Lambda_V$, then $\reip(\Lambda)=0$. Moreover $\forall \omega\in\Lambda^\complement$ we have always the convergence (\ref{cvsigma}) for all $V\in\mathcal{L}$. The collection $\mathcal{L}$ is closed under the operation of finite intersection. By Theorem 2.2 in \cite{Billingsley68} (page 14) we have $\sigma_n\Rightarrow\sigma$, $\forall \omega\in\Lambda^\complement$, which implies $\sigma_n\stackrel{a.s.}{\Rightarrow}\sigma .$ \qed

\par{\it Proof of Theorem \ref{normademc}:}\quad
We denote $c^2=b(1-b)$. In order to prove (\ref{29}) it suffices to show the following relations,
\begin{equation}\label{30}
U_n\cvdn \mathcal{N}(0,c^2),
\end{equation}

\begin{equation}\label{31}
\sqrt{n}r_m\rightarrow 0,
\end{equation}

\begin{equation}\label{32}
\displaystyle\frac{1}{n}\sum_{i=1}^n\eta_{m,i}^2-\left(\displaystyle\frac{1}{n}\sum_{i=1}^n\eta_{m,i}\right)^2\xrightarrow[N\rightarrow\infty]{P}
c^2.
\end{equation}

Since $0<\one_B(\theta_{m,i})\leq 1$, the moments $\reie|\one_B(\theta_{m,i})-\reie \one_B(\theta_{m,i})|^k$ are uniformly bounded for all $k>0$. The limit variance is 
\[c_m^2:=\reie(\eta_{m,1}-\reie\eta_{m,1})^2=\reie\eta^2_{m,1}-(\reie\eta_{m,1})^2=b+r_m-(b+r_m)^2\rightarrow c^2,\]
if we have (\ref{31}). We consider the condition of Lindeberg: for all $\varepsilon>0$,
\begin{eqnarray}\label{ingkappa}
&&\lim_{n\rightarrow\infty}\frac{1}{nc_m^2}\sum_{i=1}^n\int_{\{|\one_B(\theta_{m,i})-\reie \one_B(\theta_{m,i})|>\varepsilon\sqrt{n}c_m\}}(\one_B(\theta_{m,i})-\reie \one_B(\theta_{m,i}))^2d\reip\nonumber\\
&=&\lim_{n\rightarrow\infty}\frac{1}{c_m^2}\int_{\{|\one_B(\theta_{m,i})-\reie \one_B(\theta_{m,i})|>\varepsilon\sqrt{n}c_m\}}(\one_B(\theta_{m,i})-\reie \one_B(\theta_{m,i}))^2d\reip\nonumber\\
&\leq&\lim_{n\rightarrow\infty}\frac{1}{c_m^2}(\reie(\one_B(\theta_{m,i})-\reie \one_B(\theta_{m,i}))^4)^{1/2}(\reip\{|\varkappa_{m,1}-a_m|>\varepsilon\sqrt{n}\sigma_m\})^{1/2}\nonumber\\
&=&0.
\end{eqnarray}
The relation (\ref{30}) follows from the application of central limit theorem for the triangular array $\{\one_B(\theta_{m,i})-\reie \one_B(\theta_{m,i})\}$.

It is easy to see that (\ref{32}) follows from the application of Proposition \ref{lemgnc2} to triangular array $\{\eta_{m,i}^2\}$. It remains to establish the relation (\ref{31}). In fact, we have $|r_m|\leq C\max(m^{-1}, m^{-(\beta-\alpha)/\alpha})$ from the following lemma. Taking an arbitrarily small positive constant $\varepsilon$ and
\[n=N^{\frac{2\zeta}{1+2\zeta}-\varepsilon}, \; m=N^{\frac{1}{1+2\zeta}+\varepsilon},\]
we get $\sqrt{n}r_{m}\rightarrow 0$, the theorem is proved.\qed

\begin{lemme}\label{1}
If the condition (\ref{2ndm}) is satisfied, then
\begin{equation}\label{33}
|r_m|\leq C\max(m^{-1}, m^{-(\beta-\alpha)/\alpha}).
\end{equation}
\end{lemme}
\proof In order to prove (\ref{33}), we need to show that
\begin{equation}\label{34}
\reip\{\theta_{m,i}\in B\}=\sigma(B)+R_m,
\end{equation}
with the remainder term $R_m=O(\max(m^{-1}, m^{-(\beta-\alpha)/\alpha}))$. Let us denote
\[G_m(x)=\reip\left\{\max_{1\leq i\leq m-1}\|\xi_i\|\leq x\right\}.\]
Using the definition of $\theta_{m,i}$, it is not difficult to see that
\[\reip\{\theta_{m,i}\in B\}=m\int_0^\infty\reip\left\{\frac{\xi_1}{\|\xi_1\|}\in B, \|\xi_1\|>r\right\}G_m(dr).\]
Let $\widetilde{G}_m(x)=G_m(xm^{1/\alpha}).$ Assumption (\ref{2ndm}) implies (we recall that $\sigma(S)=1$) that for large $s$,
\begin{equation}\label{35}
\reip\{\|\xi\|>s\}=s^{-\alpha}+Cs^{-\beta}+o(s^{-\beta}).
\end{equation}
Therefore, it is easy to get the relation
\[\lim_{m\rightarrow\infty}\widetilde{G}_m(x)=G_0(x)=\left\{\begin{array}{ll}\exp({-x^{-\alpha}}),&
x>0,\\
0,& x\leq 0.\end{array}\right.\] 
Using (\ref{2ndm}) and the fact that $\int_0^\infty y^{-\alpha}dG_0(y)=1$, we have (\ref{34}) with $R_m=\sum_{i=1}^4 R_{m,i}$, where
\begin{eqnarray*}
R_{m,1}&=&m\int_0^s\reip\left\{\frac{\xi_1}{\|\xi_1\|}\in B,
\|\xi_1\|>r\right\}dG_m(r),\\
R_{m,2}&=&-\sigma(B)\int_o^{s'}y^{-\alpha}dG_0(y),\\
R_{m,3}&=&\sigma(B)\int_{s'}^\infty
y^{-\alpha}d(\widetilde{G}_m(y)-G_0(y)),\\
R_{m,4}&=&Cm^{-(\beta-\alpha)/\alpha}\int_{s'}^\infty
y^{-\beta}d\widetilde{G}_m(y).
\end{eqnarray*}
Here $s'=sm^{-1/\alpha}$ and we shall choose $s$ later. It is easy to see that
\[R_{m,1}\leq m(1-\reip\{\|\xi_1\|>s\})^m=m(1-hm^{-1})^m\leq me^{-\frac{1}{2}h},\]
where $h=h(m,s)=m\reip\{\|\xi_1\|>s\}\geq\frac{1}{2}ms^{-\alpha}$. We have used (\ref{35}) for the last inequality. Thus, if we choose
\[s=\left(\frac{m}{K\ln m}\right)^{1/\alpha},\]
with sufficiently large $K$, then we get
\begin{equation}\label{36}
R_{m,1}=o(m^{-1}).
\end{equation} 
Simple calculations show that
\begin{equation}\label{37}
R_{m,2}=o(m^{-1}).
\end{equation}
The main remainder term is $R_{m,3}$ and to estimate it we must first estimate the difference $\widetilde{G}_m(y)-G_0(y)$. A rather simple expansion of logarithmic function gives the following estimates which are sufficient for our purposes.

\begin{lemme}\label{2} (\cite{Davydov99} Lemma 2)
Let $\xi_i, i\geq 1$ be i.i.d. random vectors satisfying (\ref{35}). Then for $y>cm^{-1/\alpha}$
\begin{equation}\label{38}|\widetilde{G}_m(y)-G_0(y)|\leq
C(\alpha,\beta)\exp(-y^{-\alpha})(m^{-(\beta-\alpha)/\alpha}y^{-\beta}+m^{-1}y^{-2\alpha})
\end{equation}
and
\begin{equation}\label{39}
\sup_y|\widetilde{G}_m(y)-G_0(y)|=C(\alpha,\beta)\max(m^{-1},m^{-(\beta-\alpha)/\alpha}).
\end{equation}
\end{lemme}

Now we can estimate the term $R_{m,3}$. Integrating by parts, we get
\begin{equation}\label{40}
|R_{m,3}|=\sigma (B)(R_{m,3}^{(1)}+R_{m,3}^{(2)}),
\end{equation}
where
\[R_{m,3}^{(1)}=s'^{-\alpha}|\widetilde{G}_m(s')-G_0(s')|,\]
\[R_{m,3}^{(2)}=\alpha\int_{s'}^\infty|\widetilde{G}_m(y)-G_0(y)|y^{-\alpha-1}dy.\]
Since $s'=(K\ln m)^{-1/\alpha}>C m^{-1/\alpha}$, we can use (\ref{38}) to estimate both quantities $R_{m,3}^{(i)}, i=1,2$. After some simple calculations, we get
\[R_{m,3}^{(1)}=o(m^{-1}),\]
\[R_{m,3}^{(2)}\leq C(\alpha,\beta)\max(m^{-1},m^{-(\beta-\alpha)/\alpha}).\]

In a similar way we estimate $R_{m,4}$:
\begin{equation}\label{42}
R_{m,4}=Cm^{-(\beta-\alpha)/\alpha}\int_{s'}^\infty
y^{-\beta}d\widetilde{G}_m(y)=Cm^{-(\beta-\alpha)/\alpha}(R_{m,4}^{(1)}+R_{m,4}^{(2)}),
\end{equation}
where
\[R_{m,4}^{(1)}=\int_{s'}^\infty y^{-\beta}dG_0(y),\]
\[R_{m,4}^{(2)}=\int_s^\infty y^{-\beta}d(\widetilde{G}_m(y)-G_0(y)).\]
It is easy to see that
\begin{equation}\label{43}
R_{m,4}^{(1)}\leq C(\alpha,\beta)
\end{equation}
and $R_{m,4}^{(2)}$ can be estimated in a similar way to $R_{m,3}$:
\begin{equation}\label{44}
R_{m,4}^{(2)}\leq C(\alpha,\beta)\max(m^{-1},
m^{-(\beta-\alpha)/\alpha}).
\end{equation}
Collecting (\ref{36}), (\ref{37}), and (\ref{40})-(\ref{44}) we get (\ref{33}). 
\qed

\bibliographystyle{plain}

\end{document}